\documentclass[11pt]{article}
\title{Natural Isomorphism from a Linear Map's Image to Complement of Nullspace}
\author{H. N. Friedel}
\usepackage{amsmath}
\usepackage{latexsym}
\usepackage{graphicx}    
\begin{document}
\maketitle
\begin{abstract}
There is a natural isomorphism from image to complement of nullspace, for a bounded linear map from a real Banach space onto a closed subspace of a real Hilbert space.  This generalizes Riesz representation (self-duality of Hilbert space).  The isomorphism helps solve the pressure equation of fluid dynamics.
\end{abstract}
Notation here holds throughout.  A bound linear map $A$, from a real Banach space $X$ to a real Hilbert space $H$, has closed image $\mathcal{I} \! \mathit{m}(A) \,$.  (``Bound'' means ``bounded''.)  The complement of nullspace $\mathcal{N} (A)$, or conullspace
$\mathcal {N}^{\bot} (A) \,$, is
$\{ f \in X^* \, : \, fx=0 ~\, \mathrm{if} ~\, Ax=0 \, \}$; $X^*$ is dual-space.  $A^{t\!r}$ is transpose.  Let $\, h \in H \,$; $ \, (h | A) \,$ denotes the function $\{ X \ni x \mapsto (h | Ax) \}$ in $\, \mathcal {N}^{\bot} (A) \,$.  The promised isomorphism is $ \, \tilde{A} : \mathcal{I} \! \mathit{m}(A) \ni h \mapsto  (h | A) \in \mathcal {N}^{\bot} (A) $.\\
\textbf{Note 1.} $\tilde{A}$ is an isomorphism (of Banach spaces).  $\|\tilde{A}\| \, = \, \|A\| \,$. \\
\textbf{Proof.}  Clearly, $\tilde{A}$ is linear.  To show $\tilde{A}$ is $1\!\!:\!\!1$, we will show $\,h=0$, if $\tilde{A}h=0$.  Since $h \in \mathcal {I} \! \mathit{m}(A)$, there is $p \in X \,$ with $Ap=h$.  $\; 0 = \tilde{A}h_{\big| p} = (h|Ap) = \|h\|^2$; hence $h=0$, as promised. \\
$~~~~~$Next, we will show $\tilde{A}$ is onto.  Let $f \in \mathcal {N}^{\bot} (A)$, seek $h \in \mathcal {I} \! \mathit{m}(A)$, with $f = \tilde{A}h = (h|A)$.  Banach's closed-image theorem gives $\mathcal {N}^{\bot} (A) = \mathcal {I} \! \mathit{m}(A^{t\!r})$.  So $\,f = A^{t\!r} \, \phi$, for some $\phi \in H^*$; Riesz-representation gives $h_0 \in H$, with $(h_0 \, | \cdot) = \phi$.  Also, because $ \mathcal {I} \! \mathit{m}(A)$ is a closed subspace, it gives a direct-sum: $ \mathcal {I} \! \mathit{m}(A) \oplus  \mathcal{I} \! \mathit{m}^{\perp}(A) = H \,$.
Thus $\, h_0 \, = \, h + h_{\bot} \,$ with $\, h \in \mathcal{I} \! \mathit{m}(A) \,$, $\, h_{\bot} \in \mathcal{I} \! \mathit{m}^{\perp}(A)$.  To verify $ \, f=(h|A)$, let $x \in X$, and see:
$$ fx \; = \; (A^{t\!r} \, \phi)x \; = \; \phi (Ax) \; = \; (h_0 | Ax) \; = \; (h + h_{\bot} | Ax) \; = \; (h|Ax) \; .$$
$~~~~~$Next, we will see $\, \| \tilde{A} \| \leq \| A \| \,$.  If $h \in \mathcal {I} \! \mathit{m}(A)$, then
$$ \| \tilde{A} h \| \; = \; \| \, (h|A) \, \| \; = \; \sup_{\|x\|=1} | \, (h | Ax) \, | \; \leq \; \|h\| \, \| A \| \, , $$
because $ | \, (h|Ax) \, | \leq \|h\| \, \|Ax\| \leq \|h\| \, \|A\| \, \|x\| = \|h\| \, \|A\| \; $. \\
$~~~~~$To conclude ($\, \|\tilde{A}\| \, = \, \|A\| \,$), it remains to show $\, \|\tilde{A}\| \, \geq \, \|A\| \,$; this will come by taking supremum (over $p$) of the following bound:
\begin{equation}
\| \tilde{A} \| \; \geq \; \| Ap \| \quad \mathrm{if} \quad p \in X \; , ~~~ \|p\|=1 \; . \label{A-bound}
\end{equation}
To prove (\ref{A-bound}), let $p \in X$, $\|p\|=1$; then
\begin{equation}
\| \tilde{A} \| \, \| Ap \| \; \geq \; \| \tilde{A} (Ap) \| \; = \; \| \, (Ap|A) \, \| \; = \;  \sup_{\|x\|=1} | \, (Ap | Ax) \, | \; \geq \; | \, (Ap|Ap) \, | \; = \; \|Ap\|^2 \; . \label{A-bound2}
\end{equation}
If $\|Ap\|=0$, then certainly $\| \tilde{A} \| \; \geq \| Ap \|$; if $\|Ap\|>0$, then divide (\ref{A-bound2}) by $\|Ap\|$, to get (\ref{A-bound}).  \textbf{Done.} \\
\textit{Remarks.}  $\bullet$ Scalars must be real.  Try complex scalars, with scalar-product conjugate-linear in the first entry; then $\{ \tilde{A}: h \mapsto (h|A) \}$ is not linear. \\
$\bullet$ The Fredholm alternative $\{ \mathcal{I} \! \mathit{m}(A) \, = \, ^{\bot} \! \mathcal{N}(A^{t\!r}), \; \mathcal{I} \! \mathit{m}(A^{t\!r}) \, = \, \mathcal {N}^{\bot} (A) \}$, gives isomorphism $\{ \tilde{A} : \mathcal{I} \! \mathit{m}(A) \to \mathcal {N}^{\bot} (A) \}$ a dual expression: $\{ \tilde{A} : \, ^{\bot} \! \mathcal{N}(A^{t\!r}) \to \mathcal{I} \! \mathit{m}(A^{t\!r}) \}$, conullspace maps to image, for the transpose map.  Here, of course, $^{\bot} \! \mathcal{N}(A^{t\!r}) = \{ p \in H : (p|h)=0 ~\, \mathrm{if} ~\, (h|A)=0 \}$. \\
$\bullet$ $\tilde{A} \; = \; A^{t\!r} \, J_{\big| \mathcal {I} \! \mathit{m}(A)} \,$, $\,$where $J$ is the duality map:  $H \ni h \mapsto (h| \cdot ) \in H^* \,$. \\
$\bullet$ The isomorphism (image to conullspace) includes well-known facts.  Riesz representation ($H$ isomorphic to $H^*$) is the case where $X=H$ and $A=I$ (identity).  Also, since $ \mathcal {N}^{\bot} (A) = \mathcal {I} \! \mathit{m}(A^{t\!r}) $, the isomorphism links $ \mathcal {I} \! \mathit{m}(A)$ to $\mathcal {I} \! \mathit{m}(A^{t\!r})$: the images of a map and its transpose are naturally isomorphic; in particular, their dimensions match (row-rank equals column-rank, for matrices).  $\Box $ \\
\textbf{Acknowledgment.}  \textit{Note 1} builds on the ``generalized Riesz theorem'' of [Z], Chapter 5: if $\,f \in \mathcal {N}^{\bot} (A)$, then $\,f = (h|A)$ for some $\,h \in H$.  This follows from \textit{Note 1}, since $\,\tilde{A}$ is onto.  (More is true: $h$ may be chosen (uniquely) in $ \mathcal{I} \! \mathit{m}(A)$, yielding an isomorphism.)  $\Box $

Besides $\tilde{A}$, another natural isomorphism involving $ \mathcal {I} \! \mathit{m}(A)$ is the coset-map, $\hat{A}: X \big/ \mathcal {N}(A) \to \mathcal {I} \! \mathit{m}(A) $, $\, \hat{A}(x + \mathcal {N}(A)) = Ax \,$.  Then $\tilde{A} \hat{A}$ is a natural isomorphism from $X \big/ \mathcal {N}(A)$ to $\mathcal {N}^{\bot} (A) \,$.  $\tilde{A} \hat{A} (x + \mathcal {N}(A)) \, = \, (Ax|A) \,$. \\
\textbf{Note 2.} $\| \hat{A} \| \,=\, \|A\|$, $\| \tilde{A} \hat{A} \| \,=\, \|A\|^2 \,$. \\
\textbf{Proof.} If $\| x + \mathcal {N}(A) \| = 1$, then for each $\epsilon >0$, there is $n \in \mathcal {N}(A) \, $ with $1 \leq \| x+n \| \leq 1 + \epsilon$.  See
$$ \| \hat{A}(x + \mathcal {N}(A)) \| \; = \; \|A(x)\| \; = \; \|A(x+n)\| \; \leq \; \|A\| \|x+n\| \; \leq \; (1+\epsilon) \|A\| \; . $$
Arbitrary smallness of $\epsilon$ gives $\| \hat{A} \| \, \leq \, \| A \| \,$. \\
$~~~~~$For each $\epsilon >0$, there is $x \in X \,$ with $\|x\|=1$ and $\, 0 \leq \|A\| - \|Ax\| \leq \epsilon \,$.  See
$$ \| x + \mathcal {N}(A) \| \leq \|x\| = 1 ~ , \qquad 0 \; \leq \; \|A\| \; - \; \| \hat{A}(x + \mathcal {N}(A)) \| \; = \; \|A\| \; - \; \| Ax \| \; \leq \; \epsilon \; . $$
Arbitrary smallness of $\epsilon$ gives $\| \hat{A} \| \, = \, \| A \| \,$.  So also $\| \tilde{A} \hat{A} \| \, \leq \, \| \tilde{A} \| \, \| \hat{A} \| \, \leq \, \|A\|^2 \,$; and
$$ 0 \; \leq \; \| A \|^2 - \| \tilde{A} \hat{A} (x + \mathcal {N}(A)) \| \; = \; \|A \|^2 - \| (Ax|A) \| \; \leq \; \|A\|^2 - |(Ax|Ax)| \; = $$
$$ = \; \|A\|^2 - \|Ax\|^2 \; = \; \big( \|A\| + \|Ax\| \big) \, \big( \|A\| - \|Ax\| \big) \; \leq \; 2 \|A\| \epsilon \; . $$
Arbitrary smallness of $\epsilon$ gives $\| \tilde{A} \hat{A} \| \, = \, \| A \|^2 \, $. \textbf{Done.} \\
\textit{Remark.} The natural isomorphism, $\tilde{A} \hat{A} : X \big/ \mathcal {N}(A) \to \mathcal {N}^{\bot} (A) \,$, generalizes the following familiar isomorphism, for the case where $X=H$, with closed subspace $N$, and $A$ is the ortho-projector onto $N^{\bot}$ (with $\mathcal {N}(A)=N$).  Because here $A$ is self-adjoint with $A^2=A$, here $\tilde{A} \hat{A}$ identifies with the isomorphism $X/N \ni x+N \mapsto Ax \in N^{\bot}$.  $\Box$

The isomorphism (image to conullspace) simplifies a standard argument in fluid dynamics, as in [G] (or [Z]).  Consider flow in a region $\Omega$, a non-empty bound open connected subset of $\mathbf{R}^3$, with smooth border.  As usual, $\mathcal{L}_2 (\Omega)$, $W_2^1 (\Omega)$, $\mathring{W}_2^1 (\Omega)$ denote (respectively) the Lebesgue space of square-summable functions (modulo null measure, on $\Omega$), the Sobolev space of functions with square-summable weak-rates, and its subspace of functions vanishing on the border.
$\{ \mathcal{L}_2 \}^3$ means $\mathcal{L}_2 (\Omega) \times  \mathcal{L}_2 (\Omega) \times \mathcal{L}_2 (\Omega)$, viewed as a Hilbert space of vector-fields on $\Omega$.  $\Gamma$ denotes the subspace of gradients; $G \in \Gamma$ iff $G = -\nabla p$ for some function $p \in W_2^1$; in context of fluid dynamics, view $G$ as density of external force, and $p$ as pressure.  Given a certain force-gradient $G$, we seek its pressure $p$; to ensure uniqueness of pressure, require also $\int_{\Omega}p =0$.  The isomorphism (image to conullspace) will neatly solve the pressure equation $\{ -\nabla p =G \,, \; \int p =0 \}$, and show continuity of the solution's dependence on data.

Let $\nabla \! \cdot$ denote the divergence-map from $\{ {\mathring{W}_2^1} \}^3$ to $\mathcal{L}_2$.  Write $\mathcal{L}_2^0 \, = \, \{f \in \mathcal{L}_2 : \int f = 0  \}$.  The Divergence Theorem gives $\mathcal{I} \! \mathit{m}(\nabla \! \cdot) \subseteq \mathcal{L}_2^0 \,$; in fact, $\mathcal{I} \! \mathit{m}(\nabla \! \cdot) \, = \, \mathcal{L}_2^0 \,$, as proven in [G].  $\mathcal{L}_2^0 \,$ is closed (because it equals null-space of continuous function $\{ f \mapsto \int f \}$); hence $\nabla \! \cdot$ has closed image, and \textit{Note 1} applies (with $A = \nabla \! \cdot$).

We will see $\Gamma$ imbeds in $\mathcal{N}^{\bot} (\nabla \! \cdot)$, if we identify a gradient $G$ with the function $G^* = (G|\cdot)$.  Recall Helmholtz:  each vector-field in the complement $\Gamma ^{\bot}$ equals the limit of a sequence of smooth vector-fields with divergence zero (``di-null''), vanishing near the region's border (``border-null'').  If $V \in \mathcal{N} (\nabla \! \cdot)$, then $V$ is di-null and border-null, hence $V \in \Gamma^{\bot}$, and $0=(G|V)=G^* \, V$, $\, G^* \in \mathcal{N}^{\bot} (\nabla \! \cdot)$.

Since $\, G^* \in \mathcal{N}^{\bot} (\nabla \! \cdot) \,$, the inverse-isomorphism $\big( \tilde{A}^{-1} \big)$ of \textit{Note 1} promises just-one $p \in \mathcal{I} \! \mathit{m}(\nabla \! \cdot) \, = \, \mathcal{L}_2^0 \,$, with $G^* = (p|\nabla \! \cdot)$; the map $\{ G \mapsto p \}$ is linear, continuous, and $1\!\!:\!\!1$.  It remains only to show $p \in W_2^1$ and $-\nabla p = G$.  Write $G=(g_1 \, , \, g_2 \, , \, g_3 )$; let $\phi$ be any compactly-supported smooth function on $\Omega$; denote by $V$ the vector-field $(\phi,0,0) \in \{ {\mathring{W}_2^1} \}^3 \,$.  Then
$$ \int p \, \partial _1 \phi \; = \; (p\, \big| \, \nabla \! \cdot V) \; = \; (G|V) \; = \; \int g_1 \, \phi \; , $$
which implies $-\partial _1 p = g_1 \,$, weakly ($W_2^1$).  Likewise prove $-\partial _2 p = g_2 \,$, $-\partial _3 p = g_3 \,$; $-\nabla p = G$.

\textbf{References} \\
$\mathrm{[G]} \ \ \ $Galdi, G.$\ \ $\textit{An Introduction to the Mathematical Theory of the Navier-Stokes Equations}, vol. 1.  Springer, 1994.\\
$\mathrm{[Z]} \ \ \ $Zeidler, E.$\ \ $\textit{Applied Functional Analysis: Main Principles and Their Applications}.  Springer, 1995.
\end{document}